\documentstyle[11pt,amssymb,amstex]{amsart}

\newtheorem{thm}{Theorem}[section]

\newtheorem{cor}[thm]{Corollary}

\def\a{\alpha}
\def\o{\omega}
\def\O{\Omega}
\begin{document}

\title[Mean ergodic theorems ]
{Mean ergodic theorems in Hilbert-Kaplansky spaces}

\maketitle
\begin{center}
\author{ F.A.Shahidi $^{\rm a}$ \footnote{Corresponding author. Email: farruh.shahidi@@gmail.com} and I.G.Ganiev $^{\rm b}$ }\\\vspace{12pt}
\address{$^{\rm a}${\em{Faculty of Information and Communication Technology,
International Islamic University Malaysia, P.O Box 10, 50728, Kuala Lumpur, Malaysia}};

$^{\rm b}${\em{Faculty of Engineering, International Islamic
University Malaysia, P.O Box 10, 50728, Kuala Lumpur, Malaysia}}}

\end{center}

\begin{abstract}

We prove the mean ergodic theorem of von Neumann in a
Hilbert~---Kaplansky space. We also prove a multiparameter,
modulated, subsequential and a weighted mean ergodic theorems in a
Hilbert~---Kaplansky space \vskip 0.3cm \noindent {\it
Mathematics Subject Classification}: 37A30, 47A35, 46B42, 46E30, 46G10.
\end{abstract}

\section{Introduction}

It is well known fact that for a contracting map $T$ on a Hilbert
space $H$, the ergodic average
$\frac1n\sum\limits_{i=0}^{n-1}T^ix$ is norm convergent for all
$x\in H$ \cite{kren}, \cite{Neum}. This theorem is referred as
mean ergodic theorem or von Neumann ergodic theorem. Since then,
various generalizations of mean ergodic theorem were given. For
example, in \cite{BLRT}, the convergence of a modulated and a
subsequential ergodic averages  have been studied. Also, in
\cite{linweber}, the convergence of weighted ergodic averages in a
Hilbert space have been given. An analogue of mean ergodic theorem
in $L_2$ for multiple contractions is due to T.Tao \cite{Tao}
under assumption that transformations commute.

 Banach~---Kantorovich spaces were firstly introduced by A. V. Kantorovich
in 1938 (see, for example \cite{Kant0}, \cite{Kant} ), which have
a rich applications in analysis. Later, the theory of
Banach~---Kantorovich spaces was developed in \cite{gut0},
\cite{gut1}, \cite{kus1}, \cite{kus2}. One of the interesting
problems in a Banach~---Kantorovich space is to study the
convergence of ergodic averages. However, just few results are
known in this direction. For example, an analogue of individual
ergodic theorem for positive contractions on a
Banach~---Kantorovich lattice $L_p(\hat{\nabla}, \hat{\mu})$ has
been given in \cite{chg1}. Later in \cite{zch1}, this result was
extended to an Orlicz~---Kantorovich space. In addition,
"zero-two" law for positive contractions on Banach~---Kantorovich
lattice $L_p(\hat{\nabla}, \hat{\mu})$ has been proven in
\cite{Gm}. An analogue of Doob's martingale convergence theorem in
a Banach~---Kantorovich lattice $L_p(\hat{\nabla}, \hat{\mu})$ was
given in \cite{gan2}. So, the study of the convergence of ergodic
averages in Banach~---Kantorovich spaces is no doubt of interest.

 Since the work of Kaplansky \cite{Kapl}, a special case of
Banach~---Kantorovich space, Hilbert~---Kaplansky spaces
(AW*-modules) have been introduced. Noncommutative algebras
consummated as  subalgebras of the algebra  of bounded operators
on a Hilbert~---Kaplansky space were considered by A.G. Kusraev
\cite{kus1},\cite{kus2}. In \cite{GanArz}, a Hilbert~---Kaplansky
space over $L^0$ is represented as a measurable bundle of Hilbert
spaces $L^0$ bounded operators are given as a measurable bundle of
bounded operators in layers.

There are several approaches that may be used to get the ergodic
type theorems. One of them is a direct method. We can repeat all
the steps provided in the proof of classical Banach or Hilbert
spaces, taking into consideration the distinctions of $L^0-$
valued norms. Another way is to use boolean analysis developed in
\cite{kus1}, which gives a possibility to reduce AW*-modules to
Hilbert spaces. Finally, a Hilbert~---Kaplansky space can be
represented as a measurable bundle of classical Hilbert spaces,
which is based on the existence of respective liftings.

The first method is ineffective, because we need to repeat all the
known steps of proofs for classical Hilbert spaces, modifying
these steps into $L^0-$ valued inner product. The second method is
connected with the use of sufficiently complicated apparatus of
Boolean analysis, realization of which requires a huge preliminary
work, connected with the establishment of interrelation of
Hilbert~---Kaplansky spaces in an ordinary and boolean models of
the set theory. More effective way, in our opinion, is the use of
the third method. This is because, the theory of measurable
bundles of Banach as well as Hilbert  spaces has been developed
sufficiently well (see, for example \cite{kus2}).

In the present paper we give a description of contractions in a
Hilbert~---Kaplansky space. We prove an analogue of von Neumann
ergodic theorem and its multiparameter analogue in a
Hilbert~---Kaplansky space. Besides,we study the convergence
modulated, subsequential and weighted ergodic averages. To do so,
we use a theory of measurable bundles.

The paper is organized as follows. We give some necessary
notations and give the definition for Hilbert~---Kaplansky space
in the next section. In section 3 we give the description of
contractions and unitary operators in a Hilbert~---Kaplansky
space. We also prove the mean ergodic theorem and multiparameter
mean ergodic theorem in a Hilbert~---Kaplansky space. Finally, in
section 4 we study the convergence of modulated, subsequential and
weighted ergodic averages in a Hilbert~---Kaplansky space.

\section{Preliminaries}

Let $(\Omega,\Sigma,\mu)$ be a space with a complete finite measure,
and $L^{0}=L^{0}(\Omega)$ be an algebra of classes of  complex measurable functions in $(\Omega,\Sigma,\mu)$.

Consider a vector space $H$ over the complex numbers $\mathbb{C}$.
A transformation $||\cdot||:H\rightarrow L^0$ is said to be a
\textit{vector} or an \textit{$L^0-$ valued norm} on $H$ if it the
following conditions are fulfilled: 1) $||x||\ge 0$ for all $x\in
H$; $||x||=0\Leftrightarrow x=0;$ 2) $||\lambda
x||=|\lambda|||x||$ for all $\lambda\in C$ and $x\in H$; 3)
$||x+y||\le ||x||+||y||$ for all $x,y\in H.$

A pair $(H, ||\cdot||)$ is said to be \textit{lattice-normed}
space over $L^0.$ Lattice-normed space $H$ is said to be
\textit{disjunctively decomposable}, or $d-$ decomposable if the
following condition holds: for any $x\in H$ and disjunctive
elements $e_1,e_2\in L^0$ satisfying $||x||=e_1+e_2,$ there exist
$x_1,x_2\in H$ such that $x=x_1+x_2$ with $||x_1||=e_1$ and
$||x_2||=e_2.$ A net $x_{\a}$ in $H$ is said to be a
\textit{$(bo)-$ convergent} to an element $x\in H$ if there exists
a decreasing net $(e_{\gamma})_{\gamma\in\Gamma}$ with
$inf_{\gamma\in\Gamma}e_{\gamma}=0$ such that for any
$\gamma\in\Gamma$ there is $\a=\a(\gamma)$ such that
$||x-x_{\a}||\le e_{\gamma}$ for all $\a\ge\a(\gamma).$ A
lattice-normed space is said to be \textit{$(bo)-$ complete} if
any fundamental net ${x_{\a}}$ in H is $(bo)-$ convergent to an
element of $H.$ Any $d$-decomposable $(bo)-$ complete
lattice-normed space is said to be a
\textit{Banach~---Kantorovich} space.

\textbf{Definition 1} \cite{kus2}(\cite{kus1}). A transformation
$\langle\cdot,\cdot\rangle:H\times H\rightarrow L^{0}$ is called
an \textit{$L^{0}$-valued inner product}, if for any $x,y,z\in H$
and $\alpha\in \mathbb{C}$, the following conditions hold :

         1) $\langle x,y\rangle\geqslant 0,\langle x,x\rangle=0\Leftrightarrow x=0$;
         2) $\langle x,y+z\rangle=\langle x,y\rangle+\langle x,z\rangle$;
         3) $\langle\alpha x,y\rangle=\alpha\langle x,y\rangle$;
           4) $\langle x,y\rangle=\overline{\langle y,x\rangle}$.

It is known \cite{kus2}, that  $||x||=\sqrt{\langle x,x\rangle}$
defines an $L^{0}$-valued norm in $H$. If $(H,||\cdot||)$ is a
Banach~---Kantorovich space, then
 $(H,\langle\cdot,\cdot\rangle)$ is said to be a \textit{Hilbert~---Kaplansky} space.
  Examples of such spaces can be found in \cite{kus1}, \cite{kus2}.

 Let ${\cal H}$ be a map that assigns some Hilbert space $\cal H(\omega)$
 to any point $\omega \in \Omega.$ Function $u$, defined a.e. in $\Omega$, with the values
 $ u(\omega) \in \cal H(\omega)$, for all $\o$ in the domain  $  dom (u) $  of $u$, is said to be
 a \textit{section} on $\cal H.$ For the set of sections $L$, following \cite{gut1}, we call a pair
 $({\cal H} ,L)$ a \textit{measurable bundle}  of Hilbert spaces over $\Omega$, if

1) $\lambda_1 c_1+\lambda_2 c_2 \in L$ for all $\lambda_1,\lambda_2\in C$
and $c_1,c_2\in L$, where $\lambda_1 c_1+\lambda_2 c_2 :\omega\in \textmd  {dom} $ $(c_1)$ $\bigcap \textmd  {dom} $ $(c_2)$ $\rightarrow\lambda_1 c_1(\omega)+\lambda_2 c_2(\omega)$;

2) a function $||c||: \omega \in \textmd  {dom} $ $(c)$
$\rightarrow||c(\omega)||_{\cal H(\omega)}$ is measurable for all
$c\in L$;

3) for any $\omega\in\Omega$ the set $\{c(\omega):c\in L,\omega\in \textmd  {dom} $ $(c)$\}
  is dense in $\cal H(\omega)$.

A section $s$ is called a \textit{step-section}, if
$s(\omega)=\sum\limits_{i=1}^n \chi_{A_i}(\omega)c_i(\omega)$,
where $c_i\in L,A_i\in \Sigma,i=\overline{1,n}.$ A section $u$ is
called \textit{measurable}, if there exists a sequence
$(s_n)_{n\in \mathbb{N}}$ of step-sections such that
$||s_n(\omega)-u(\omega)||_{\cal H(\omega)}\rightarrow 0$ for
almost all $\omega\in\Omega$.

Let $M(\Omega,{\cal H})$ be the set of all measurable sections. By
$L^{0}(\Omega,{\cal H})$ we denote the factor space of
$M(\Omega,{\cal  H})$ with respect to a.e. equality. By
$\widehat{u}$ we denote the class from $L^{0}(\Omega,{\cal H})$,
containing the section $u.$ Note that the function $\omega
\rightarrow ||u(\omega)||_{\cal H(\omega)}$ is measurable for any
$u\in M(\Omega,{\cal H})$, and therefore, the function
$(u(\o),v(\o))_{\cal H(\o)}=
 \frac{1}{4}(||u(\o)+v(\o)||^{2}_{\cal
 H(\o)}-||u(\o)-v(\o)||^{2}_{\cal H(\o)})$ is also measurable for all $u,v\in
 M(\Omega,{\cal H}).$

 We denote by $\langle\widehat{u},\widehat{v}\rangle$ the element
of $L^{0}$ containing  $(u(\o),v(\o))_{\cal H(\o)}$. Clearly,
$\langle\cdot,\cdot\rangle$ is an $L^{0}$-valued inner product. We
also denote by $||\widehat{u}||_{{\cal H}(\o)}$ the element of
$L^{0}$ containing the function $||u(\o)||,$ for any $u\in
M(\Omega,{\cal H})$. Then
$||\widehat{u}||^{2}=\widehat{||u(\o)||_{\cal
H(\o)}^{2}}=\widehat{(u(\o),u(\o))}_{\cal
H(\o)}=\langle\widehat{u},\widehat{u}\rangle$, that is
     $||u||=\sqrt{\langle\widehat{u},\widehat{u}\rangle}$.

From theorem 4.1.14  (\cite{gut1} page 144) $(L^{0}(\Omega,{\cal
H}),||\cdot||)$ is a Banach~---Kantorovich space. That is why
$(L^{0}(\Omega,{\cal
 H}),\langle\cdot,\cdot\rangle)$ is a Hilbert~---Kaplansky space over $L^{0}$.

\textbf{Definition 2.} The collection $\{T_{\o}:{\cal
H(\o)}\rightarrow {\cal H(\o)},\ \o\in\O \}$ of linear operators
is called a \textit{measurable bundle of linear operators}, if
$T_{\o}x(\o)\in M(\O, {\cal H})$ for any $x\in M(\O, {\cal H}).$

$H$  be a Hilbert~---Kaplansky space. An operator $T:H\rightarrow
H$ on a Hilbert~---Kaplansky space $H$  is said to be an
\textit{$L^0-$ linear}, if $T(\alpha x+\beta y)=\alpha T(x)+\beta
T(y)$ for all $\alpha, \beta\in L^0(\O)$ and $x,y\in H.$ $L^0-$
linear operator is called an \textit{$L^0-$ bounded}, if for any
bounded set $B$ in $L^0,$ the set $T(B)$ is bounded as well. For
the $L^0-$ bounded operator $T$ we define its norm as
$||T||=sup\{||Tx||:\ ||x||\le\textbf{ 1}\},$ where $\textbf{1}$ is
an identity element of $L^0.$

If every $T_{\o}$ is contraction (unitary), then $\{T_{\o}:\ \o\in\O \}$ is called a measurable bundle of
contractions (unitary operators).

\begin{thm}\cite{GanKud} Let $\{T_{\o}:{\cal H(\o)}\rightarrow {\cal H(\o)},\ \o\in\O \}$ be a measurable bundle of linear operators.
Then the operator $\hat{T}: L^0(\O, {\cal H})\rightarrow L^0(\O, {\cal H}),$ defined by $\hat{T}\hat{x}=\widehat{T_{\o}x(\o)}$
is $L^0-$ linear and $L^0-$ bounded operator.
\end{thm}

Let ${\cal L}^{\infty}(\O)$ be the set of all bounded measurable
functions on $\O$ and $L^{\infty}(\O)$ be the factor space of
${\cal L}^{\infty}(\O)$  with respect to a.e. equality. By ${\cal
L}^{\infty}(\Omega, {\cal H})$ we denote the set of those points
$||u||\in M(\O, {\cal H})$ for which $||u(\o)||_{H(\o)}\in {\cal
L}^{\infty}(\Omega)$ and by $L^{\infty}(\O, {\cal H})$ we denote
the factor space of ${\cal L}^{\infty}(\O, {\cal H})$ with respect
to equality a.e.

Let $p:L^{\infty}(\O) \rightarrow {\cal L}^{\infty}(\O)$ be a
lifting (see \cite{gut1},\cite{kus2}).

\textbf{Definition 3.} A map $l:L^{\infty}(\Omega, {\cal H})
\rightarrow {\cal L}^{\infty}(\Omega, {\cal H}) $ is said to be a
\textit{vector valued lifting} associated with  lifting $p$ if for
all $\hat{u},\ \hat{v}\in L^{\infty}(\Omega, {\cal H}) $ the
following conditions hold:

1) $l(\hat{u})\in\hat{u}, \ \ dom(u)=\Omega;$

2) $||l(\hat{u})(\o)||_{H(\o)}=p(||\hat{u}||)(\o);$

3)$l(\hat{u}+\hat{v})=l(\hat{u})+l(\hat{v});$

4) $l(\lambda\hat{u})=p(\lambda)l(\hat{u})$ for $\lambda\in L^{\infty}(\Omega, {\cal H});$

5) $l(\hat{u}\hat{v})=l(\hat{u})l(\hat{v});$

6) $\{l(\hat{u})(\o): \hat{u}\in L^{\infty}(\Omega,
\mathcal{H})\}$ is dense in $\mathcal{H(\o)},$ for all $\o\in\O.$

\begin{thm}\cite{GanArz} For any Hilbert~---Kaplansky space $H$ over $L^0$ there exists
a measurable bundle of Hilbert spaces
 $({\cal H}, L)$, with vector valued lifting, such that
$H$ isometrically isomorphic to $L^0(\Omega, {\cal H})$.

\end{thm}

\section{Mean ergodic theorem in a Hilbert~---Kaplansky space}

In this section we  give a description of contractions  and prove an
analogue of von Neumann ergodic theorem in a Hilbert~---Kaplansky
space. Besides, we prove the convergence multiparameter  ergodic
averages.

An $L^0-$ linear and $L^0-$ bounded operator $T:L^0(\Omega, {\cal
H})\rightarrow L^0(\Omega, {\cal H})$ is called a
\textit{contraction (unitary) in $L^0(\O, \mathcal{H})$} if
$||T||\le \textbf{1} \ (||T||= \textbf{1})$ The following theorem
describes contractions in Hilbert~---Kaplansky space.

\begin{thm} For any contraction  $T:L^0(\Omega, {\cal H})\rightarrow L^0(\Omega, {\cal H})$
and for all $\o\in\Omega$ there exists a measurable bundle of
contractions $T_{\o}:\mathcal{H(\o)}\rightarrow \mathcal{H(\o)}$
such that
$$(T\hat{u})(\o)=T_{\o}(u(\o))$$
for almost all $\o\in\O$ and for all $\hat{u}\in L^0(\Omega,
\mathcal{H}).$
\end{thm}

\begin{pf}
Since $T$ is a contraction in $L^0(\O, {\cal H})$, then
$||T\hat{u}||\le ||\hat{u}||.$ Let $\hat{u}\in L^{\infty}(\Omega,
{\cal H}).$ Then we have $||\hat{u}||\in L^{\infty}(\O)$, hence
$||T\hat{u}||\in L^{\infty}(\Omega),$ which means that
$T\hat{u}\in L^{\infty}(\Omega, {\cal H}).$

Let $T_{\o}(l(\hat{u})(\o))=l(T\hat{u}(\omega)),$ where $l$ is a
vector valued lifting on $L^{\infty}(\Omega, \mathcal{H}),$
associated with a lifting $p$. From

$$||T_{\o}(l(\hat{u})(\o))||_{{\cal H(\o)}}=||l(T\hat{u})(\o)||_{{\cal H(\o)}}=p(||T\hat{u}||)(\o)\le$$$$
\le p(||\hat{u}||)(w)=||l(\hat{u})(\o)||_{{\cal H(\o)}}$$ we imply
that the linear operator

$$T_{\o}: \{l(\hat{u})(\o): \hat{u}\in L^{\infty}(\Omega, \mathcal{H})\}\rightarrow \mathcal{H(\o)}$$
is well defined  and bounded. Moreover $||T_{\o}||\le 1.$

Since $\{l(\hat{u})(\o): \hat{u}\in L^{\infty}(\Omega, {\cal
H})\}$ is dense in ${\cal H(\o)},$ then for all $\o\in\Omega,$
$T_{\o}$ can be extended to operator $T_{\o}:{\cal
H(\o)}\rightarrow {\cal H(\o)}$ preserving the norm. That is why
we get contractions $T_{\o}:{\cal H(\o)}\rightarrow {\cal H(\o)}$
for all $\o\in\Omega.$

According to definition of $T_{\omega}$ and property 1) of a vector valued lifting we have
$$T_{\o}(l(\hat{u})(\o))=l(T\hat{u})(\o)=(T\hat{u})(\o)$$
 for almost all $\o\in\Omega.$

Now, let $\hat{u}\in L^0(\Omega, {\cal H}).$ Since
$L^{\infty}(\Omega, {\cal H})$ is (bo)-dense in $L^0(\Omega, {\cal
H})$ (see \cite{gan1}), then there is a sequence $\hat{u}_n\in
L^{\infty}(\Omega, {\cal H})$ such that
$||\hat{u}_n-\hat{u}||\rightarrow 0.$ Then
$||\hat{u}_n(\o)-\hat{u}(\o)||_{{\cal H(\o)}}\rightarrow 0$ for
almost all $\o\in\O.$ From
$$T(\hat{u})=\lim\limits_{n\to\infty}T(\hat{u}_n)$$
we get
$$||T_{\o}(\hat{u}_n(\o))-T(\hat{u})(\o)||_{{\cal H(\o)}}= ||T_{\o}(\hat{u}_n)(\o)-T(\hat{u})(\o)||_{{\cal H(\o)}}\rightarrow 0$$
for almost all $\o\in\O.$ Therefore,
$T(\hat{u})(\o)=\lim\limits_{n\to\infty}T_{\o}(u_n(\o))$ for
almost all $\o\in\O.$ On the other hand, the continuity of
$T_{\o}$ yields
$\lim\limits_{n\to\infty}T_{\o}(u_n(\o))=T_{\o}(u(\o)).$ Hence for
every $\hat{u}\in L^0(\Omega, {\cal H})$ we have
$T(\hat{u})(\o)=T_{\o}(u(\o))$ for almost all $\o\in\O.$

\end{pf}
\begin{cor}For any  unitary operator $U:L^0(\Omega, {\cal H})\rightarrow L^0(\Omega, {\cal H})$
and for all $\o\in\Omega$ there exists a measurable bundle of
unitary operators $U_{\o}:{\cal H(\o)}\rightarrow {\cal H(\o)}$
such that
$$(U\hat{u})(\o)=U_{\o}(u(\o))$$
for almost all $\o\in\Omega$ and for all $\hat{u}\in L^0(\Omega,
\mathcal{H}).$
\end{cor}

For $\hat{u}\in L^0(\O, \mathcal{H})$ and $T:L^0(\Omega,
{\cal H})\rightarrow L_0(\Omega, {\cal H})$ we define $A_n(T,\
\hat{u})=\frac1n\sum\limits_{k=0}^{n-1}T^k\hat{u}.$

The following theorem is an analogue of von Neumann ergodic
theorem for Hilbert~---Kaplansky space.

\begin{thm} If $T: L^0(\O, {\cal H})\rightarrow L^0(\O, {\cal H})$ is
a contraction, then  $A_n\hat{u}\rightarrow \hat{P}(\hat{u})$ in
$L^0(\Omega, {\cal H})$ for all $\hat{u}\in L^0(\Omega, {\cal
H})$, where $\hat{P}$ is a projection in $L^0(\O, {\cal H})$.
\end{thm}

\begin{pf} Let $T$ be a contraction in $L^0(\O, {\cal H})$. According to Theorem 3.1 there exist the corresponding contractions $T_{\o}$ in ${\cal H(\o)}.$
Then we have
$$A_n(T,\ \hat{u})(\o)=(\frac1n\sum\limits_{k=0}^{n-1}T^k\hat{u})(\o)=\frac1n\sum\limits_{k=0}^{n-1}T_{\o}^k(u(\o))=A_n(u(\o))$$
for any $\hat{u}\in L^0(\Omega, {\cal H})$ and almost all
$\o\in\O.$

From mean ergodic theorem (see \cite{kren}) there exists a
projection $P(\o)$  such that $A_n(u(\o))\rightarrow P(\o)(u(\o))$
in ${\cal H(\o)}.$ Since $A_n(u(\o))$ is a measurable section,
then from Gutman's  theorem \cite{gut1}, it follows that the
section $P(\o)(u(\o))$ is measurable. Further, since projections
$P(\o)$ map a measurable bundle to a measurable bundle, then
Theorem 2.1 implies the existence of $L^0-$ linear and $L^0-$
bounded operator $\hat{P}: L^0(\O, {\cal H})\rightarrow L^0(\O,
{\cal H})$ given by
$$\hat{P}(\hat{u})=\widehat{P(\o)(u(\o))}.$$

Finally, from

$$||A_n(T,\ \hat{u})-\hat{P}(\hat{u})||=||\widehat{A_n(u(\o))-P(\o)(u(\o))}||_{{\cal H(\o)}}\rightarrow 0$$
we get $A_n(T,\ \hat{u})\rightarrow \hat{P}(\hat{u})$ in
$L^0(\Omega,{\cal H}).$ We show that $\hat{P}$ is a projection.
Indeed, since $P(\o)$ is a projection, then we have
$P(\o)^2=P(\o)$ and $P(\o)^*=P(\o),$ where $P(\o)^*$ is the
adjoint of $P(\o).$ Since $P(\o)$ is  measurable bundle of
operators, then  $\{P(\o)^*: \o\in\O \}$ and $ \{P(\o)^2: \o\in\O
\}$ are measurable bundles of operators. We have
$$\hat{P}^2(\hat{u})=\widehat{P(\o)^2u(\o)}=\widehat{P(\o)u(\o)}=\hat{P}(\hat{u}).$$
Similarly,
$$\hat{P}^*(\hat{u})=\widehat{P(\o)^*u(\o)}=\widehat{P(\o)u(\o)}=\hat{P}(\hat{u}).$$
So, the ergodic average tends to some projection.

\end{pf}
In the following examples we illustrate the above proved theorem.

\textbf{Example 1.} Let $\O=\mathbb{N},$ then $L^0(\O)=s-$ the
space of sequences. We define
$$ H=s[\ell_2]=\{f=(f_1, f_2,\cdots, f_n,\cdots): f_i\in\ell_2, \ \forall i\in\mathbb{N}\}.$$
 The norm is defined as $||f||=(||f_1||_{\ell_2}, ||f_2||_{\ell_2},\cdots, ||f_n||_{\ell_2}, \cdots)\in s,$
 and the inner product as $\langle f,g\rangle=((f_1,g_1), (f_2,g_2), \cdots, (f_n,g_n),\cdots)\in
 s,$ where $(\cdot, \cdot)$ is an inner product in $\ell_2.$
Then $H$ is a Hilbert~---Kaplansky space.

Now, let $f_i=(f_i^{(1)}, f_i^{(2)},\cdots,
f_i^{(n)},\cdots,)\in\ell_2$ and $Tf_i=(0,
f_i^{(1)},f_i^{(2)},\cdots).$ Define $s-$ linear operator
$\textbf{T}: s[\ell_2]\rightarrow s[\ell_2]$ by
$$\textbf{T}f=(Tf_1, Tf_2,\cdots, Tf_n,\cdots).$$
One can see that
$$||\textbf{T}f||=(||Tf_1||_{\ell_2}, ||Tf_2||_{\ell_2},\cdots)=(||f_1||_{\ell_2}, ||f_2||_{\ell_2}, \cdots)=||f||,$$
hence $\textbf{T}$ is a contraction. Then according to Theorem 3.3
$\lim\limits_{n\to\infty}\frac1n\sum\limits_{k=0}^{n-1}\textbf{T}^kf$
is norm convergent. Let us find its limit. Since
$\frac1n\sum\limits_{k=0}^{n-1}T^kf_i\rightarrow 0$ as
$n\to\infty$ in $\ell_2$ for all $i\in\mathbb{N},$ then
$$\frac1n\sum\limits_{k=0}^{n-1}\textbf{T}^kf=(\frac1n\sum\limits_{k=0}^{n-1}T^kf_1, \frac1n\sum\limits_{k=0}^{n-1}T^kf_2,\cdots)\rightarrow(0,0,\cdots)$$
as $n\to\infty.$ Therefore,
$\lim\limits_{n\to\infty}\frac1n\sum\limits_{k=0}^{n-1}\textbf{T}^kf=0$
in $s[\ell_2].$

\textbf{Example 2.} Let $([0,\ \infty), {\cal B}, m)$ be
measurable space and $L_2=L_2([0,\ \infty), {\cal B}, m)$ be a
Banach space of square integrable functions. Let $L^0[L_2]$ be the
space of all measurable function $K$ on $\O\times\ [0,\ \infty)$
containing the class of equivalency of functions $y\rightarrow
K(\o, y)$ belonging to $L_2$ such that $\o\rightarrow
||K(\o,\dot)||_{L_2}$ is measurable.

We define the norm of $K=K(\o, y)$ as $||K||=\widehat{||K(\o,
\cdot)||_{L_2}}\in L^0.$ Then $L^0[L_2]$ is a
Banach~---Kantorovich space \cite{kus2}. We define the inner
product as $\langle K_1,K_2\rangle=\widehat{(K_1(\o, y),K_2(\o,
y))_{L_2}},$ where $(\cdot ,\cdot)_{L_2}$ denotes the inner
product in $L_2.$ Then $||K||=\sqrt{\langle K,K\rangle}$ implies
that $L^0[L_2]$ is a Hilbert~---Kaplansky space.

Now let $P(t,x, B)$ be a Markov process with invariant measure
$m.$ Then the operator
$$\textbf{T}(K)=\int\limits_0^{\infty}K(\o, y)P(1,x,dy)$$
is a contraction \cite{Sadd}. Then according to Theorem 3.3
$\lim\limits_{n\to\infty}\frac1n\sum\limits_{k=0}^{n-1}\textbf{T}^k(K)$
is norm convergent in $L^0[L_2]$.

In the following example we assume $\O$ to be a countable.

\textbf{Example 3.} Let $\O$ be a countable set. Then
$L^0[L_2]=s[L_2].$ For $K\in s[L_2]$ we  define
$||K||=(||K_1||_{L_2}, ||K_2||_{L_2},\cdots
||K_n||_{L_2},\cdots)\in s$ and $\langle K,
M\rangle=((K_1,M_1)_{L_2}, (K_2,M_2)_{L_2},\cdots).$ Then $s[L_2]$
is a Hilbert~---Kaplansky space.

 We define
$\textbf{T}:s[L_2]\rightarrow s[L_2]$ as
$$\textbf{T}(K)=(T(K_1),T(K_2),\cdots)$$
where $T(K_i)=K_i(y+1).$ Since $T$ is a contraction in $L_2,$ then
$\textbf{T}$ is a  contraction in $s[L_2].$ Then
$$A_n(\textbf{T}, K)=(A_n(T, K_1), A_n(T, K_2),\cdots).$$
Since $A_n(T, K_i)$ converges as $n\to\infty$  in $L_2$ for all
$i$ then $A_n(\textbf{T}, K)$ converges in $s[L_2].$

Let $T_1,\ T_2,\cdots, T_d$ be $d$ linear operators on
$L^0(\Omega, {\cal H}).$ We define the \textit{multiparameter
ergodic average} by
$$A_{n_1,n_2,\cdots, n_d}(T_1,\ T_2,\cdots, T_d,\ \hat{u})=\frac1{n_1n_2\cdots n_d}\sum\limits_{i_1=0}^{n_1-1}\sum\limits_{i_2=0}^{n_2-1}
\cdots\sum\limits_{i_d=0}^{n_d-1}T_1^{i_1}T_2^{i_2}\cdots
T_d^{i_d}(\hat{u}).$$

\begin{thm} Let $T_1,\ T_2,\cdots, T_d$ be contractions  in $L^0(\O, {\cal H})$.
Then for any $\hat{u}\in L^0(\Omega, {\cal H}),$ the average
$A_{n_1,n_2,\cdots,n_k}(T_1,\ T_2,\cdots, T_d,\ \hat{u})$ is norm
convergent as $n_1,n_2,\cdots, n_d\rightarrow\infty$
independently.
\end{thm}

\begin{pf} We proceed by induction, noting that the theorem is true for $d=1$ by the previous theorem.
Assume that multiparameter ergodic average is norm convergent for
any $d-1$ contractions in $L^0(\O, {\cal H})$. Let
$P_1,P_2,\cdots, P_d$ be projections such that
$$\lim\limits_{n_i\to\infty}A_{n_i}(T_i,\ \hat{u})=P_i(\hat{u})\ \ i=\overline{1,d}.$$
Then $$\lim\limits_{n_1,n_2,\cdots,n_d\to\infty}A_n(T_1,\ T_2,\cdots, T_d,\ \hat{u})=P_dP_{d-1}\cdots P_1(\hat{u}).$$

Indeed,

$$||A_{n_1,n_2,\cdots, n_d}(T_1,\ T_2,\cdots, T_d,\ \hat{u})-P_dP_{d-1}\cdots P_1(\hat{u})||\le$$$$\le
||A_{n_2,\cdots, n_d}(T_2,\cdots, T_d,\
\hat{u})(A_{n_1}(T_1)-P_1)\hat{u}||+||(A_{n_2,\cdots,
n_d}(T_2,\cdots, T_d,\ \hat{u})-P_d\cdots
P_2)P_1(\hat{u})||\le$$$$ \le
||A_{n_1}(T_1)(\hat{u})-P_1(\hat{u})||+||(A_{n_2,\cdots,
n_d}(T_2,\cdots, T_d,\ \hat{u})-P_d\cdots P_2)P_1(\hat{u})||$$

Due to the assumption of induction both expressions
$||A_{n_1}(T_1)(\hat{u})-P_1(\hat{u})||$ and

$||(A_{n_2,\cdots, n_d}(T_2,\cdots, T_d,\ \hat{u})-P_d\cdots
P_2)P_1(\hat{u})||$ converge to $0$ as
$n_1,n_2,\cdots,n_d\to\infty.$ Therefore, multiparameter ergodic
average is norm convergent.

\end{pf}

\section{Modulated, subsequential and weighted ergodic theorems}

In this section we study a modulated, subsequential and a weighted
ergodic theorems in a Hilbert~---Kaplansky space.

Let $\{a_k\}_{k\ge 0}$ be a sequence of complex numbers. For
$\hat{u}\in L^0(\O, {\cal H})$ and $T:L^0(\Omega, {\cal
H})\rightarrow L^0(\Omega, {\cal H})$ we define $A_n(a_k,\ T,\
\hat{u})=\frac1n\sum\limits_{k=0}^{n-1}a_kT^k\hat{u}.$ Following
\cite{BLRT}, we call this average a \textit{modulated}.

\begin{thm} Let $\{a_k\}_{k\ge 0}$ be a sequence of complex numbers satisfying the following conditions:

a) For every complex $\lambda$ with $|\lambda|=1$ there exists $c(\lambda)$ such that

$$\frac1n\sum\limits_{k=0}^{n-1}a_k\overline{\lambda}^k\rightarrow c(\lambda)$$ as $n\to\infty$

b) $$sup_{n\ge 1}sup_{|\lambda|=1}|\frac1n\sum\limits_{k=0}^{n-1}a_k\lambda^k|<\infty.$$

Then for any contraction $T$ in  $L^0(\O, {\cal H}),$ the average
$A_n(a_k,\ T,\ \hat{u})$ converges in norm for any $\hat{u}\in
L^0(\Omega, {\cal H});$

\end{thm}

\begin{pf}  Let $T$ be a contraction in $L^0(\O, {\cal H})$ and $T_{\o}$ be the corresponding contractions in ${\cal H(\o)}.$
Then
$$A_n(a_k,\ T,\  \hat{u})(\o)=(\frac1n\sum\limits_{k=0}^{n-1}a_kT^k\hat{u})(\o)=\frac1n\sum\limits_{k=0}^{n-1}a_kT_{\o}^k(u(\o))=A_n(a_k,\ T_{\o},\ u(\o))$$
for any $\hat{u}\in L^0(\Omega, {\cal H})$ and almost all
$\o\in\O.$

Since the sequence $\{a_k\}_{k\ge 0}$ satisfies a) and b), then
according to Corollary 2.3. from \cite{BLRT}, there exists
$u^*(\o)$ such that $A_n(a_k,\ T_{\o},\ u(\o))=
\frac1n\sum\limits_{k=0}^{n-1}a_kT_{\o}^k(u(\o))\rightarrow
u^*(\o)$ in ${\cal H(\o)}$ for almost all $\o\in\O$ as
$n\to\infty.$ Note that $A_n(a_k,\ T_{\o},\ u(\o))$ is a
measurable section, and therefore according to \cite{gut1},
$u^*(\o)$ is a measurable section and $\hat{u}^*=\widehat{u^*}\in
L^0(\Omega, {\cal H}).$ Therefore, from

 $$||A_n(a_k,\ T,\ \hat{u})-\hat{u}^*||=||\widehat{A_n(a_k,\ T_{\o},\ u(\o))-u^*(\o)}||_{{\cal H(\o)}}\rightarrow 0$$
we get $A_n(a_k,\ T,\ \hat{u})\rightarrow \hat{u}^*$ in
$L^0(\Omega, {\cal H}).$

\end{pf}

\begin{cor} If the conditions of above Theorem 4.1 hold, then for any unitary operator $U$ in $L^0(\O, {\cal H})$,
 the average $A_n(a_k,\ U,\ \hat{u})$ converges in norm for any $\hat{u}\in L^0(\Omega, {\cal H}).$
\end{cor}

\begin{pf} Applying Corollary 2.3 from \cite{BLRT} and Corollary 3.2 and the
arguments given in the proof of we get the desired result.
\end{pf}

Now, we turn our attention to the study of subsequential ergodic
averages. For a sequence $\{k_j\}$, $\hat{u}\in L^0(\O, {\cal H})$
and $T:L^0(\Omega, {\cal H})\rightarrow L^0(\Omega, {\cal H})$ we
define a subsequential ergodic average $A_n(k_j,\ T,\
\hat{u})=\frac1n\sum\limits_{k=0}^{n-1}T^{k_j}\hat{u}.$

\begin{thm} Let $\{k_j\}$ be a strictly increasing sequence of positive integers satisfying
$\lim\limits_{n\to\infty}\sum\limits_{j=1}^n\lambda^{k_j}=0$ for
every $|\lambda|=1,\ \lambda\neq 1.$ Then for a given contraction
$T$ in $L^0(\O, {\cal H}),$ the average $A_n(k_j,\ T,\ \hat{u})$
converges in norm for all $\hat{u}\in L^0(\Omega, {\cal H}).$
\end{thm}

\begin{pf} Let $T$ be an $L^0-$ contraction and $T_{\o}$ be the corresponding measurable bundle of contractions in ${\cal H(\o)}.$ Then
$$A_n(k_j,\ T,\  \hat{u})(\o)=(\frac1n\sum\limits_{j=1}^{n}T^{k_j}\hat{u})(\o)=\frac1n\sum\limits_{k=0}^{n-1}T_{\o}^{k_j}(u(\o))=A_n(k_j,\ T_{\o},\ u(\o))$$
for any $\hat{u}\in L^0(\Omega, {\cal H})$ and almost all
$\o\in\O.$

The sequence $\{k_j\}$ satisfies the condition of the Theorem,
then Proposition 3.2 from \cite{BLRT} implies that there exists
$u^*(\o),$ such that $A_n(k_j,\ T_{\o},\ u(\o))=
\frac1n\sum\limits_{j=1}^{n}T_{\o}^{k_j}(u(\o))\rightarrow
u^*(\o)$ in ${\cal H(\o)}$ as $n\to\infty.$ Note that $A_n(k_j,\
T_{\o},\ u(\o))$ is a measurable section, hence $u^*(\o)$ is
measurable and $\hat{u}^*=\widehat{u^*}\in L^0(\Omega, {\cal H}).$
Therefore, from

 $$||A_n(k_j,\ T,\ \hat{u})-\hat{u}^*||=||\widehat{A_n(k_j,\ T_{\o},\ u(\o))-u^*(\o)}||_{{\cal H(\o)}}\rightarrow 0$$
we get $A_n(a_k,\ T,\ \hat{u})\rightarrow \hat{u}^*$ in
$L^0(\Omega, {\cal H}).$
\end{pf}

Let $\{w_k\}_{k\ge 0}$ be a non-null sequence of nonnegative
numbers and denote its partial sums by $W_n.$ We also define
$A_n(w_k,\ T,\
\hat{u})=\frac1{W_n}\sum\limits_{k=0}^{n-1}w_kT^k\hat{u}.$ The
following theorem is a weighted ergodic theorem in a
Hilbert~---Kaplansky space.

\begin{thm} If for any complex $\lambda$ with $|\lambda|=1$ we have
$$\frac1{W_n}\sum\limits_{k=0}^{n-1}w_k\overline{\lambda}^k\rightarrow c(\lambda)$$ as $n\to\infty.$
Then

a). For any contraction $T$ in  $L^0(\O, {\cal H}),$ the average
$A_n(w_k,\ T,\ \hat{u})$ converges in norm for any $\hat{u}\in
L^0(\Omega, {\cal H});$

b). For any contraction $T$ in  $L^0(\O, {\cal H})$
$$||A_n(w_k,\ T,\ \hat{u})-A_n( T,\ \hat{u})||\rightarrow 0.$$
\end{thm}

\begin{pf} a). This part can be proven by providing all arguments given in the proof of the previous
theorem and applying Theorem 2.1 from \cite{linweber}.

b). Let $T$ be a contraction in $L^0(\O, {\cal H})$ and $T_{\o}$
be the corresponding contractions in ${\cal H(\o)}.$ Then
$$A_n(w_k,\ T,\  \hat{u})(\o)=(\frac1{W_n}\sum\limits_{k=0}^{n-1}w_kT^k\hat{u})(\o)=\frac1{W_n}\sum\limits_{k=0}^{n-1}w_kT_{\o}^k(u(\o))=A_n(w_k,\ T_{\o},\ u(\o))$$
for any $\hat{u}\in L^0(\Omega, {\cal H})$ and almost all
$\o\in\O.$

From Corollary 2.2 of \cite{linweber} we get
$$||A_n(w_k,\ T_{\o},\  u(\o))-A_n(T_{\o},\  u(\o))||\rightarrow 0$$
in ${\cal H(\o)}$ for almost all $\o\in\O.$ Using the technique of
the previous theorems we obtain
 $$||A_n(w_k,\ T,\ \hat{u})-A_n(T,\ \hat{u})||=\widehat{||A_n(w_k,\ T_{\o},\ u(\o))-A_n(T_{\o},\ u(\o))||_{{\cal H(\o)}}}\rightarrow 0.$$

\end{pf}
Let $k_j$ be a strictly increasing sequence of positive integers;
put $w_{k_j} = 1$, and $w_k = 0$ if $k\notin\{k_j\}.$ Then the
weighted averages become the averages along the subsequence $k_j.$
More examples of weights $w_k$, satisfying the conditions of
Theorem 4.4 can be found in \cite{linweber}.


\begin{thebibliography}{99}




\bibitem{BLRT} Berend D., Lin M., Rosenblatt J., Tempelman A. Modulated and subsequential ergodic theorems
in Hilbert and Banach spaces. Ergod. Th.\& Dynam. Sys. 22(2002), 1653-1665.


\bibitem{chg1} Chilin V.I. ,Ganiev I.G. Individual ergodic theorem for contraction in Banach-Kontorovich lattice
 $L^p(\hat{\nabla},\hat{\mu}).$  Russian Mathematics,
 7(2000), 81--83.


 \bibitem{gan2} Ganiev I.G., Martingale convergence on Banach-Kontorovich lattices $L^p(\hat{\nabla},\hat{\mu}).$
 Uzbek Math. Journal, 1(2000), 18--25.


\bibitem{gan1} Ganiev I.G. Measurable bundles of lattices and their applications
Reseach on Functional Analysia and its applications, Moskov Nauka,
2006. 10--49 (Russian).

\bibitem{GanArz} Ganiev I.G., Arziev A.D. Spectrum of self-adjoint
operators in Hilber~---Kaplansky modules over $L^0.$ Proceedings
of Int. Conf.: Research on analysis, mathematical modeling and
informatics. Vladikavkaz, Institute of applied mathematics and
informatics VNS RAN, 2007, 7-23.

\bibitem{Gm} Ganiev I.G., Mukhamedov F.,  On the "Zero-Two" law for positive contractions in the Banach-Kantorovich lattice $Lp(\nabla ,\mu
)$. Comment. Math. Univ. Carolinae.  47(2006), 427--436.

\bibitem{GanKud} Ganiev I.G., Kudaybergenov K.K.  Banach theorem
on inverse operator in Banach~---Kantorovich space. Vladikavkaz
Math. J. 3(6)(2004), 321-325.


\bibitem{gut0} Gutman A.E., Banach bundles in the theory of lattice normed space, II,
Siberian Adv. Math. 3(4)(1993) 8-40.


\bibitem{gut1} Gutman A.E., Banach bundles in lattice normed spaces,
Linear operators, associated with order, Novosibirsk: Izd. IM SO RAN, 1995, 63-211.

\bibitem{Kant0} Kantorovich L. V., "On a class of functional equations," Dokl. Akad. Nauk
SSSR, 4, No. 5(1936), 211–216.
\bibitem{Kant} Kantorovich L.V., Vulih B.Z., Pinsker A.G. Functional analysis in partially ordered spaces,
Gosudarstv. Izdat. Tehn.-Teor. Lit., Moscow-Leningrad, 1950
(Russian).

\bibitem{Kapl} Kaplansky I., Modules over operator algebras, Amer. J. Math. 75(4)(1953), 839-858.


\bibitem{kren} Krengel U., Ergodic Theorems. Walter de
Grugwer.--Berlin, New-York. 1985, 357~p.


\bibitem{kus1} Kusraev A.G., Kutateladze S.S. An introduction to boolean valued analysis. Moscow: Nauka, 1985 (Russian).

\bibitem{kus2} Kusraev A.G. Dominated operators, Kluwer Acad. Publ., 2000, 446~p.



\bibitem{linweber} Lin M., Weber M., Weighted ergodic theorems and strong laws of large nambers. Ergod. Th.
and Dynam. Sys. 27(2007), 511-–543.


\bibitem{Neum} von Neumann J., Proof of quasi ergodic hypothesis,
Proc Nac.Acad. Sci. 18(1932) 70-82.

\bibitem{Sadd} Sadaddinova S.S., On semigroups of operators on the
space of Bochner measurable functions. Vestnik NUUz-Tashkent
3(2010), 169-172.


\bibitem{Tao} Tao T. Norm convergence of multiple ergodic averages for commuting transformation, Ergodic Theory
and Dynamical Systems 28(2008), 657-688.




\bibitem{zch1} Zakirov B.S. Chilin V.I., Ergodic theorems for contractions in Orlicz-Kontorovich lattices
, Siberian Math. Journal, 6(50)(2009), 1305-1318.


\end{thebibliography}
\end{document}